# CONSIDERATIONS ON NEW FUNCTIONS IN NUMBER THEORY

by Florentin Smarandache, Ph. D.
University of New Mexico
Gallup, NM 87301, USA

**Abstract**:
New functions are introduced in number theory, and for each one a general description, examples, connections, and references are given.

**Introduction.**
In this paper a small survey is presented on eighteen new functions and four new sequences, such as: Inferior/Superior f-Part, Fractional f-Part, Complementary function with respect with another function, S-Multiplicative, Primitive Function, Double Factorial Function, S-Prime and S-Coprime Functions, Smallest Power Function.

1) Let f: Z ---> Z be a strictly increasing function and x an element in R. Then:

**a) Inferior f-Part of x**,
**-------------------------------**
ISf(x) is the smallest k such that f(k) <= x < f(k+1).

**b) Superior f-Part of x**,
**-------------------------------**
SSf(x) is the smallest k such that f(k) < x <= f(k+1).

Particular cases:

*a) Inferior Prime Part*:
For any positive real number n one defines ISp(n) as the largest prime number less than or equal to n.
The first values of this function are (Smarandache[6] and Sloane[5]):
2,3,3,5,5,7,7,7,7,11,11,13,13,13,13,17,17,19,19,19,19,23,23.

*b) Superior Prime Part*:
For any positive real number n one defines SSp(n) as the smallest prime number greater than or equal to n.
The first values of this function are (Smarandache[6] and Sloane[5]):
2,2,2,3,5,5,7,7,11,11,11,11,13,13,17,17,17,17,19,19,23,23,23.

*c) Inferior Square Part*:
For any positive real number n one defines ISs(n) as the largest square less than or equal to n.

The first values of this function are (Smarandache[6] and Sloane[5]):
0,1,1,1,4,4,4,4,4,9,9,9,9,9,9,9,16,16,16,16,16,16,16,16,16,25,25.

*b) Superior Square Part*:
For any positive real number n one defines SSs(n) as the smallest square greater than or equal to n.
The first values of this function are (Smarandache[6] and Sloane[5]):
0,1,4,4,4,9,9,9,9,9,16,16,16,16,16,16,16,25,25,25,25,25,25,25,25,25,36.

*d) Inferior Cubic Part*:
For any positive real number n one defines ISc(n) as the largest cube less than or equal to n.
The first values of this function are (Smarandache[6] and Sloane[5]):
0,1,1,1,1,1,1,1,8,8,8,8,8,8,8,8,8,8,8,8,8,8,8,8,8,8,27,27,27,27.

*e) Superior Cube Part*:
For any positive real number n one defines SSs(n) as the smallest cube greater than or equal to n.
The first values of this function are (Smarandache[6] and Sloane[5]):
0,1,8,8,8,8,8,8,8,27,27,27,27,27,27,27,27,27,27,27,27,27,27,27,27,27.

*f) Inferior Factorial Part*:
For any positive real number n one defines ISf(n) as the largest factorial less than or equal to n.
The first values of this function are (Smarandache[6] and Sloane[5]):
1,2,2,2,2,6,6,6,6,6,6,6,6,6,6,6,6,6,6,6,6,6,6,24,24,24,24,24,24,24.

*g) Superior Factorial Part*:
For any positive real number n one defines SSf(n) as the smallest factorial greater than or equal to n.
The first values of this function are (Smarandache[6] and Sloane[5]):
1,2,6,6,6,6,24,24,24,24,24,24,24,24,24,24,24,24,24,24,24,24,24,24,120.

Remark 1: This is a generalization of the inferior/superior integer part of a number (floor function).

2) Let f: Z ---> Z be a strictly increasing function and x an element in R. Then:

**Fractional f-Part of x**,

-------------------------

FSf(x) = x - ISf(x),
where ISf(x) is the Inferior f-Part of x defined above.

Particular cases:

*a) Fractional Prime Part*:
FSp(x) = x - ISp(x),
where ISp(x) is the Inferior Prime Part defined above.
Example: FSp(12.501) = 12.501 - 11 = 1.501.

*b) Fractional Square Part*:
FSs(x) = x - ISs(x),
where ISs(x) is the Inferior Square Part defined above.
Example: FSs(12.501) = 12.501 - 9 = 3.501.

*c) Fractional Cubic Part*:
FSc(x) = x - ISc(x),
where ISc(x) is the Inferior Cubic Part defined above.
Example: FSc(12.501) = 12.501 - 8 = 4.501.

*d) Fractional Factorial Part*:
FSf(x) = x - ISf(x),
where ISf(x) is the Inferior Factorial Part defined above.
Example: FSf(12.501) = 12.501 - 6 = 6.501.

Remark 2.1: This is a generalization of the fractional part of a number.
Remark 2.2: In a similar way one defines:
- the Inferior Fractional f-Part:
  IFSf(x) = x - ISf(x) = FSf(x);
- and the Superior Fractional f-Part:
  SFSf(x) = SSf(x) - x;
  for example: Superior Fractional Cubic Part of 12.501
  = 27 - 12.501 = 14.499.

3) Let g: A ---> A be a strictly increasing function, and let "~" be a given internal law on A. Then we say that
**f: A ---> A is complementary with respect to the**
**---------------------------------------------**

**function g and the internal law "~"** if:
**----------------------------------**

f(x) is the smallest k such that there exists a z in A so that x~k = g(z).

Particular cases:

*a) Square Complementary Function*:
f: N ---> N, f(x) = the smallest k such that xk is a perfect square.
The first values of this function are (Smarandache[6] and Sloane[5]):
1,2,3,1,5,6,7,2,1,10,11,3,14,15,1,17,2,19,5,21,22,23,6,1,26,3,7.

*b) Cubic Complementary Function*:
f: N ---> N, f(x) = the smallest k such that xk is a perfect cube.
The first values of this function are (Smarandache[6] and Sloane[5]):
1,4,9,2,25,36,49,1,3,100,121,18,169,196,225,4,289,12,361,50.

More generally:

*c) m-power Complementary Function*:
f: N ---> N, f(x) = the smallest k such that xk is a perfect m-power.

*d) Prime Complementary Function*:
f: N ---> N, f(x) = the smallest k such that x+k is a prime.
The first values of this function are (Smarandache[6] and Sloane[5]):
1,0,0,1,0,1,0,3,2,1,0,1,0,3,2,1,0,1,0,3,2,1,0,5,4,3,2,1,0,1,0,5.

**4) S-Multiplicative Function**:

A function f : $N^*$ --> $N^*$ which,
for any (a, b) = 1, verifies f(ab) = max {f(a), f(b)};
(i.e. it reflects the main property of the Smarandache function[8]).

**5) Smarandache-Kurepa Function**:

For p prime, SK(p) is the smallest integer such that !SK(p) is divisible by p, where !SK(p) = 0! + 1! + 2! + ... + (p-1)!

For example:

| p | 2 | 3 | 7 | 11 | 17 | 19 | 23 | 31 | 37 | 41 | 61 | 71 | 73 | 89 |
|---|---|---|---|---|---|---|---|---|---|---|---|---|---|---|
| SK(p) | 2 | 4 | 6 | 6 | 5 | 7 | 7 | 12 | 22 | 16 | 55 | 54 | 42 | 24 |

**6) Smarandache-Wagstaff Function**:

For p prime, SW(p) is the smallest integer such that W(SW(p)) is divisible by p, where W(p) = 1! + 2! + ... + (p)!

For example:

| p | 3 | 11 | 17 | 23 | 29 | 37 | 41 | 43 | 53 | 67 | 73 | 79 | 97 |
|---|---|---|---|---|---|---|---|---|---|---|---|---|---|
| SW(p) | 2 | 4 | 5 | 12 | 19 | 24 | 32 | 19 | 20 | 20 | 7 | 57 | 6 |

**7) Smarandache Ceil Functions of n-th Order**:

Sk(n) is the smallest integer for which n divides Sk(n)^k.

For example, for k=2, we have:

| n | 1 | 2 | 3 | 4 | 5 | 6 | 7 | 8 | 9 | 10 | 11 | 12 | 13 | 14 | 15 | 16 |
|---|---|---|---|---|---|---|---|---|---|---|---|---|---|---|---|---|
| S2(n) | 2 | 4 | 3 | 6 | 10 | 12 | 5 | 9 | 14 | 8 | 6 | 20 | 22 | 15 | 12 | 7 |

**8) Pseudo-Smarandache Function**:

Z(n) is the smallest integer such that 1 + 2 + ... + Z(n) is divisible by n.

For example:

| n | 1 | 2 | 3 | 4 | 5 | 6 | 7 |
|---|---|---|---|---|---|---|---|
| Z(n) | 1 | 3 | 2 | 3 | 4 | 3 | 6 |

**9) Smarandache Near-To-Primordial Function**:

* * *

SNTP(n) is the smalest prime such that either p - 1, p , or p + 1 is divisible by n,

where $p^*$, of a prime number p, is the product of all primes less than or equal to p.

For example:

| n | 1 | 2 | 3 | 4 | 5 | 6 | 7 | 8 | 9 | 10 | 11 | ... | 59 ... |
|---|---|---|---|---|---|---|---|---|---|---|---|---|---|
| SNTP(n) | 2 | 2 | 2 | 5 | 3 | 3 | 3 | 5 | ? | 5 | 11 | ... | 13 ... |

## 10) Double-Factorial Function:

SDF(n) is the smallest number such that SDF(n)!! is divisible by n, where the double factorial
    m!! = 1x3x5x...xm, if m is odd;
and m!! = 2x4x6x...xm, if m is even.

For example:

| n | 1 | 2 | 3 | 4 | 5 | 6 | 7 | 8 | 9 | 10 | 11 | 12 | 13 | 14 | 15 | 16 |
|---|---|---|---|---|---|---|---|---|---|---|---|---|---|---|---|---|
| SDF(n) | 1 | 2 | 3 | 4 | 5 | 6 | 7 | 4 | 9 | 10 | 11 | 6 | 13 | 14 | 5 | 6 |

## 11) Primitive Functions:

Let p be a positive prime.

$S_p$ : N ---> N, having the property that $(S_p(n))!$ is divisible by $p^n$,

and it is the smallest integer with this property.

For example:

$S_3(4)$ = 9, because 9! is divisible by 3^4, and it is the smallest one with this property.

These functions help computing the Smarandache Function.

**12) Smarandache Function**:

$S : N \to N$, $S(n)$ is the smallest integer such that $S(n)!$ is divisible by n.

**13) Smarandache Functions of the First Kind**:

$S_n : N^* \to N^*$

i) If $n = u^r$ (with u = 1, or u = p prime number), then

$S_n(a) = k$, where k is the smallest positive integer such that k! is a multiple of $u^{ra}$;

ii) If $n = p_1^{r_1} \cdot p_2^{r_2} \cdots p_t^{r_t}$, then $S_n(a) = \max_{1<=j<=t} \{ S_{p_j^{r_j}}(a) \}$.

**14) Smarandache Functions of the Second Kind**:

$S^k : N^* \to N^*$, $S^k(n) = S_n(k)$ for k in $N^*$,

where $S_n$ are the Smarandache functions of the first kind.

**15) Smarandache Function of the Third Kind**:

$S_a^b(n) = S_{a_n}(b_n)$, where $S_{a_n}$ is the Smarandache function of the

first kind, and the sequences $(a_n)$ and $(b_n)$ are different from the following situations:

i) $a_n = 1$ and $b_n = n$, for $n$ in $N^*$;

ii) $a_n = n$ and $b_n = 1$, for $n$ in $N^*$.

**16) S. Prime Functions** are defined as follows:

$P : N \to \{0, 1\}$, with

$$P(n) = \begin{cases} 0, & \text{if } n \text{ is prime;} \\ 1, & \text{otherwise.} \end{cases}$$

For example P(2) = P(3) = P(5) = P(7) = P(11) = 0, whereas P(0) = P(1) = P(4) = P(6) = ... = 1.

More general:

$P_k : N^k \to \{0, 1\}$, where k is an integer >= 2, and

$$P_k(n_1, n_2, \ldots, n_k) = \begin{cases} 0, & \text{if } n_1, n_2, \ldots, n_k \text{ are all prime numbers;} \\ 1, & \text{otherwise.} \end{cases}$$

**17) S. Coprime Functions** are similarly defined:

$C_k : N^k \to \{0, 1\}$, where k is an integer >= 2, and

$$C_k(n_1, n_2, \ldots, n_k) = \begin{cases} 0, & \text{if } n_1, n_2, \ldots, n_k \text{ are coprime numbers;} \\ 1, & \text{otherwise.} \end{cases}$$

**18) The Smallest Power Function**:

SP(n) is the smallest number m such that $m^k$ is divisible by n, where $k \geq 2$ is given.

The following sequence SP(n) is generated:
1, 2, 3, 2, 5, 6, 7, 4, 3, 10, 11, 6, 13, 14, 15, 4, 17, 6, 19, 10, 21, 22, 23, 6, 5, 26, 3, 14, 29, 30, 31, 4, 33, 34, 35, 6, 37, 38, 39, 20, 41, 42, ...

Remarks:

If p is prime, then SP(p) = p.
If r is square free, then SP(r) = r.

If $n = (p_1^{s_1}) \times \dots \times (p_k^{s_k})$ and all $s_i \leq p_i$, then SP(n) = n.

If $n = p^s$, where p is prime, then:

$$SP(n) = \begin{cases} p, & \text{if } 1 \leq s \leq p; \\ p^2, & \text{if } p+1 \leq s \leq 2p^2; \\ p^3, & \text{if } 2p^2+1 \leq s \leq 3p^3; \\ \dots\dots\dots & \\ p^t, & \text{if } (t-1)p^{(t-1)}+1 \leq s \leq tp^t . \end{cases}$$

Generally, if $n = (p_1^{s_1}) \times \dots \times (p_k^{s_k})$, with all $p_i$ prime, then:

$SP(n) = (p_1^{t_1}) \times \dots \times (p_k^{t_k})$, where

$t_i = u_i$ if $(u_i - 1)p_i^{(u_i - 1)}+1 \leq s_i \leq u_i p_i^{u_i}$

for $1 \leq i \leq k$.

Particular cases:

*a) A second function (k=2)*:

1, 2, 3, 2, 5, 6, 7, 4. 3, 10, 11, 6, 13, 14, 15, 4, 17, 6, 19, 10, 21, 22, 23, 12, 5, 26, 9, 14, 29, 30, 31, 8, 33, ...

( S2(n) is the smallest integer m such that $m^2$ is divisible by n )

*b) A third function (k=3)*:

1, 2, 3, 2, 5, 6, 7, 8, 3, 10, 11, 6, 13, 14, 15, 4, 17, 6, 19, 10, 21, 22, 23, 6, 5, 26, 3, 14, 29, 30, 31, 4, 33, ...

( S3(n) is the smallest integer m such that $m^3$ is divisible by n )

**19) A 3n-digital subsequence**:

13, 26, 39, 412, 515, 618, 721, 824, 927, 1030, 1133, 1236, ...
(numbers that can be partitioned into two groups such that the second is three times biger than the first)

**20) A 4n-digital subsequence**:

14, 28, 312, 416, 520, 624, 728, 832, 936, 1040, 1144, 1248, ...
(numbers that can be partitioned into two groups such that the second is four times biger than the first)

**21) A 5n-digital subsequence**:

15, 210, 315, 420, 525, 630, 735, 840, 945, 1050, 1155, 1260, ...
(numbers that can be partitioned into two groups such that the second is five times bigger than the first)

**22) Sequences of Sub-sequences**

For all of the following sequences:

*a) Crescendo Sub-sequences*:

1, 1, 2, 1, 2, 3, 1, 2, 3, 4, 1, 2, 3, 4, 5, 1, 2, 3, 4, 5, 6,

1, 2, 3, 4, 5, 6, 7, 1, 2, 3, 4, 5, 6, 7, 8, . . .

*b) Decrescendo Sub-sequences*:

1, 2, 1, 3, 2, 1, 4, 3, 2, 1, 5, 4, 3, 2, 1, 6, 5, 4, 3, 2, 1,
7, 6, 5, 4, 3, 2, 1, 8, 7, 6, 5, 4, 3, 2, 1, . . .

*c) Crescendo Pyramidal Sub-sequences*:

1,   1, 2, 1,   1, 2, 3, 2, 1,   1, 2, 3, 4, 3, 2, 1,
1, 2, 3, 4, 5, 4, 3, 2, 1,   1, 2, 3, 4, 5, 6, 5, 4, 3, 2, 1,   . . .

*d) Decrescendo Pyramidal Sub-sequences*:

1,   2, 1, 2,   3, 2, 1, 2, 3,   4, 3, 2, 1, 2, 3, 4,
5, 4, 3, 2, 1, 2, 3, 4, 5,   6, 5, 4, 3, 2, 1, 2, 3, 4, 5, 6,   . . .

*e) Crescendo Symmetric Sub-sequences*:

1, 1,   1, 2, 2, 1,   1, 2, 3, 3, 2, 1,   1, 2, 3, 4, 4, 3, 2, 1,
1, 2, 3, 4, 5, 5, 4, 3, 2, 1,   1, 2, 3, 4, 5, 6, 6, 5, 4, 3, 2, 1, . . .

*f) Decrescendo Symmetric Sub-sequences*:

1, 1,   2, 1, 1, 2,   3, 2, 1, 1, 2, 3,   4, 3, 2, 1, 1, 2, 3, 4,
5, 4, 3, 2, 1, 1, 2, 3, 4, 5,   6, 5, 4, 3, 2, 1, 1, 2, 3, 4, 5, 6,  . . .

*g) Permutation Sub-sequences*:

1, 2,   1, 3, 4, 2,   1, 3, 5, 6, 4, 2,   1, 3, 5, 7, 8, 6, 4, 2,
1, 3, 5, 7, 9, 10, 8, 6, 4, 2,   1, 3, 5, 7, 9, 10, 8, 6, 4, 2, . . .

Find a formula for the general term of the sequence.

Solutions:

For purposes of notation in all problems, let

$$a(n)$$

denote the n-th term in the complete sequence and

$$b(n)$$

the n-th subsequence.  Therefore, a(n) will be a number and b(n) a sub-sequence.

a) Clearly, b(n) contains n terms. Using a well-known summation formula, at the end of b(n) there would be a total of

$$\frac{n(n + 1)}{2}$$

terms.

Therefore, since the last number of b(n) is n,

$$a((n(n+1))/2) = n.$$

Finally, since this would be the terminal number in the sub-sequence

$$b(n) = 1, 2, 3, . . . , n$$

the general formula is

$$a(((n(n+1)/2) - i) = n - i$$

for n >= 1 and 0 <= i <= n - i.

b) With modifications for decreasing rather than increasing, the proof is essentially the same. The final formula is

$$a(((n(n+1))/2) - i) = 1 + i$$

for n >= 1 and 0 <= i <= n - 1.

c) Clearly, b(n) has 2n - 1 terms. Using the well-known formula of summation

$$1 + 3 + 5 + . . . + (2n - 1) = n^2.$$

the last term of b(n) is in position $n^2$ and $a(n^2) = 1$. The largest number in b(n) is n, so counting back n - 1 positions, they increase in value by one each step until n is reached.

$$a(n^2 - i) = 1 + i, \quad \text{for } 0 <= i <= n-1.$$

After the maximum value at n-1 positions back from $n^2$, the values decrease by one. So at the nth position back, the value is n-1, at the (n-1)st position back the value is n-2 and so forth.

$$a(n^2 - n - i) = n - i - 1$$

for 0 <= i <= n - 2.

d) Using similar reasoning

$$a(n^2) = n \quad \text{for } n >= 1$$

and

$$a(n^2 - i) = n - i, \quad \text{for } 0 <= i <= n-1$$

$$a(n^2 - n - i) = 2 + i, \text{ for } 0 <= i <= n-2.$$

e) Clearly, b(n) contains 2n terms. Applying another well-known summation formula

$$2 + 4 + 6 + . . . + 2n = n(n+1), \text{ for } n >= 1.$$

Therefore, a(n(n+1)) = 1. Counting backwards n-1 positions, each term decreases by 1 up to a maximum of n.

a((n(n+1))-i) = 1 + i, for 0 <= i <= n-1

The value n positions down is also n and then the terms decrease by one back down to one.

a((n(n+1))-n-i) = n - i, for 0 <= i <= n - 1.

f) The number of terms in b(n) is the same as that for (e). The only difference is that now the direction of increase/decrease is reversed.

a((n(n+1))-i) = n - i, for 0 <= i <= n-1.

a((n(n+1))-n-i) = 1 + i, for 0 <= i <= n - 1.

g) Given the following circular permutation on the first n integers.

$$\varphi_n = \begin{pmatrix} 1 & 2 & 3 & 4 & . & . & . & n-2 & n-1 & n \\ 1 & 3 & 5 & 7 & . & . & . & 6 & 4 & 2 \end{pmatrix}$$

Once again, b(n) has 2n terms. Therefore,

a(n(n+1)) = 2.

Counting backwards n-1 positions, each term is two larger than the successor

a((n(n+1))-i) = 2 + 2i,  for 0 <= i <= n-1.

The next position down is one less than the previous and after that, each term is again two less the successor.

a((n(n+1))-n-i) = 2n - 1 - 2i, for 0 <= i <= n-1.

As a single formula using the permutation

$a((n(n+1)-i) = \varphi_n(2n-i)$, for 0 <= i <= 2n-1.

## GENERAL REFERENCES

*BOOKS*:

*ABSTRACTS*:

*PROPOSED PROBLEMS AND SOLUTIONS*:

Atanassov, K., "B l cenoiu, Ion; Seleacu, Vasile / Some Properties of Smarandache Functions of the Type I", in <Zentralblatt fuer Mathematik>, Berlin, Germany, 831 - 45, 11007, 1996/06.

Atanassov, K., "Ashbacher, Charles / Some Problems on Smarandache Function", in <Zentralblatt fuer Mathematik>, Berlin, Germany, 831 - 45, 11008, 1996/06.

Atanassov, K., "B l cenoiu, Ion; Popescu, Marcela; Seleacu, Vasile / About the Smarandache Square's Complementary Function", in <Zentralblatt fuer Mathematik>, Berlin, Germany, 831 - 45, 11009, 1996/06.

Atanassov, K., "Tomi , Tiberiu Florin / Some Remarks Concerning the Distribution of the Smarandache Function", in <Zentralblatt fuer Mathematik>, Berlin, Germany, 831 - 45, 11010, 1996/06.

Atanassov, K., "R descu, E.; R descu, N.; Dumitrescu, C. / Some Elementary Algebraic considerations Inspired by the Smarandache Function", in <Zentralblatt fuer Mathematik>, Berlin, Germany, 831 - 46, 11011, 1996/06.

Atanassov, K., "B l cenoiu, Ion; Dumitrescu, Constantin / Smarandache Functions of the Second Kind", in <Zentralblatt fuer Mathematik>, Berlin, Germany, 831 - 46, 11012, 1996/06.

Atanassov, K., "Popescu, Marcela; Popescu, Paul / The Problem of Lipschitz Condition", in <Zentralblatt fuer Mathematik>, Berlin, Germany, 831 - 46, 11013, 1996/06.

Atanassov, K., "Dumitrescu, Constantin / A Brief History of the 'Smarandache Function'", in <Zentralblatt fuer Mathematik>, Berlin, Germany, 831 - 46, 11014, 1996/06.

Atanassov, K., "Dumitrescu, C.; Seleacu, V. / Some Notions and Questions in Number Theory", in <Zentralblatt fuer Mathematik>, Berlin, Germany, 840 - 70, 11001, 1996/15.

Atanassov, K., "Smarandache, Florentin / Only Problems, Not Solutions! 4th ed", in <Zentralblatt fuer Mathematik>, Berlin, Germany, 840 - 70, 11002, 1996/15.

*CITATIONS*:

Gardner, Martin, "Lucky Numbers and 2187", in <The Mathematical Inteligencer>, USA, Vol. 19, No. 2, 1997, p. 29.

----, Mathematics Calendar, "First International Conference on Smarandache Type Notions in Number Theory", in <Notices of the American Mathematical Society>, Providence, RI, USA, September 1996, p. 1057.

----, "Smarandache Notions", in <Ulrich's International Periodicals Directory 1997>, Mathematics, USA, pp. 4396, 8074, 9117.

----, "Smarandache function", in <Subject Headings>, Library of Congress, 17th Edition, Vol. III (K-P), Washington, D.C., 1994, p. 3530.

*INTERNATIONAL CONGRESS*:

The First International Conference on Smarandache Type Notions in Number Theory, August 21-24, Department of Mathematics, University of Craiova, Romania; This Conference has been organized by Dr. C. Dumitrescu & Dr. V. Seleacu, under the auspices of UNESCO.